\documentclass[12pt]{amsart}
\usepackage{times}
\usepackage{amsmath}
\usepackage{amsthm}
\usepackage{tikz}
\usetikzlibrary{matrix,arrows}

\newcommand{\bA}{\mathbb{A}}
\newcommand{\bone}{1_{Y}}
\newcommand{\sB}{\mathcal{B}}
\newcommand{\sH}{\mathcal{H}}
\newcommand{\sHo}{\mathcal{H}_{0}}
\newcommand{\piU}{(\pi, U)}
\newcommand{\Linf}{L^{\infty}}
\newcommand{\sK}{\mathcal{K}}
\newcommand{\LY}{\Linf (Y, \nu)}
\newcommand{\cent}{\mathbf{Z}(\pi(A)')}

\newcommand{\gspace}{(\Gamma, \mu)}
\newcommand{\cstar}{C^*}
\newcommand{\cp}{(A, G, \sigma)}

\theoremstyle{plain}
\newtheorem{thm}{Theorem}

\newtheorem{lemma}[thm]{Lemma}
\newtheorem{cor}[thm]{Corollary}

\theoremstyle{remark}
\newtheorem{remark}[thm]{Remark}

\begin{document}
\title[Covariant Representations involving Compact Groups] {Covariant Representations of $\cstar$-algebras and Their Compact Automorphism Groups}
\author{Firuz Kamalov}
\address{Mathematics Department, University of Nebraska-Lincoln, Lincoln, NE, 68588-0130}
\email{s-fkamalo1@math.unl.edu}
\date{\today}

\begin{abstract}
Let $G$ be a compact group. Let $(\Gamma, \mu)$ be a standard
Borel $G$-measure space. We show that the group action on
 $(\Gamma, \mu)$ is transitive if and only if it is ergodic. Using this result,
we show that every irreducible covariant representation of a
$\cstar$-dynamical system $\cp$ is induced from a stability
group. In addition, we show that $\cp$ satisfies strong-EHI.

\end{abstract}
\maketitle

\section{Introduction}
Let $G$ be a locally compact group, $A$ a $\cstar$-algebra, and
$\sigma$ a point-wise norm continuous homomorphism of $G$ into the
automorphism group of $A$ then we call the triple $\cp$  a
$\cstar$-dynamical system. Given a $\cstar$-dynamical system we can
construct the crossed product $\cstar$-algebra $A \times_{\sigma} G$
that encodes the action of $G$ on $A$. It is well known that there
exists a one to one correspondence between the set all covariant
representations of the system $\cp$ and the set of all
*-representations of $A\times_{\sigma} G$. Therefore, the study of
representations of $A\times_{\sigma} G$ is equivalent to that of
covariant representations of $\cp$.

Our goal is to study induced covariant representations of
systems involving compact groups. The study of induced
representations was initiated by Mackey in \cite{M1, M2} in the
context of unitary representations of locally compact groups.
Using Mackey's approach Takesaki extended the theory to crossed
products in \cite{T}. Subsequently, Rieffel recast that theory
in terms of Hilbert modules and Morita equivalence with
\cite{R}. It follows from Proposition 5.4 in \cite{W1} that the
construction of induced representations for crossed products by
Rieffel is equivalent to that of Takesaki.

The importance of induced representations arises from the fact that
the fundamental structure of a crossed product $A\times_{\sigma} G$
is reflected in the structure of the orbit space for the $G$-action
on Prim $A$ together with the subsystems $(A, G_{P}, \sigma)$ where
$G_{P}$ is the stability group at $P\in$ Prim $A$. In particular,
one gets a complete description of the primitive ideal space and its
topology for transformation group $\cstar$-algebra
$C_{0}(X)\times_{lt} G$ when $G$ is abelian. In many important cases
we also get a characterization of when $A\times_{\sigma} G$ is GCR
or CCR. Williams presents all these results and more in his book
\cite{W1}.

Although induced representations have been studied extensively there
remains a considerable gap in the theory. We outline below two
questions for which answers are not known. Using structure theorems
obtained in this paper we give a positive answer to both questions
in the case of separable $\cstar$-dynamical systems with compact
groups.

One of the key ingredients in building the connection between
Prim $A\times_{\sigma} G$ and the $G$-action on Prim $A$ is
establishing that every primitive ideal of $A\times_{\sigma} G$
is induced from a stability group (\cite{W1}; p.235). The
latter result was conjectured by Effros and Hahn, and systems
for which the conjecture holds are called EH-regular. The proof
that the Effros-Hahn conjecture holds is due Gootman, Rosenberg
and Sauvageot and it is one of the major results in the theory
(see Chapters 8 and 9 in \cite{W1} for the proof of the GRS
theorem and its applications). There exists a stronger notion
of EH-regularity namely the requirement that every irreducible
representation of $A\times_{\sigma} G$ is induced from a
stability group. The latter requirement is known to hold for
many dynamical systems (\cite{W1}; Theorem 8.16) but the
general case, to our knowledge, remains open.

Another natural question that arises when studying induced
representations of a system $\cp$ is when an irreducible
representation of a subsystem $(A, G_{P}, \sigma)$ induces to an
irreducible representation of $\cp$. Following the nomenclature
proposed by Echterhoff and Williams in \cite{W2}, we say that $\cp$
satisfies strong Effros-Hahn Induction Property (strong-EHI), if,
for each primitive ideal $P$ of $A$ and a covariant irreducible
representation $(\pi, U)$ of $(A, G_{P},\sigma)$ such that
ker($\pi$) = $P$ the corresponding induced representation of $\cp$
is irreducible. A very nice summary of the results regarding the
(strong)-EHI property can be found in \cite{W2}.

In this paper we use Takesaki's approach to the theory of
induced representations for crossed products. As in \cite{T} we
will often assume basic countability conditions although most
of the results in Section 3 do hold in greater generality. If
$G$ is a second countable, locally compact group acting on a
separable $\cstar$-algebra $A$ then we call $\cp$ a separable
system.

In Section 2, we give the background about topological and Borel
dynamical systems necessary for Section 3. In Section 3, we study
Borel dynamical systems. In particular, we prove that if $G$ is a
compact group and $\gspace$ is an ergodic standard Borel $G$-measure
space then $G$ acts transitively on $\gspace$. Note that the last
statement is not true in general. For instance, the action of
$\mathbb{Z}$ on $\mathbb{T}$ by an irrational rotation is ergodic
but it is not transitive.

In Section 4, we study covariant representation $\piU$ of a
system $\cp$ on a Hilbert space $\sH$. Given a covariant
representation $\piU$ and a system of imprimitivity $\bA$ for
$\piU$ there exists an essentially unique standard Borel
$G$-measure space $\gspace$ such that $\Linf \gspace$ is
isomorphic to $\bA$. If $G$ acts ergodically on $\bA$ then the
corresponding action on $\gspace$ is also ergodic (\cite{M};
Theorem 3). In particular, by the result mentioned in the
previous paragraph $G$ acts transitively on $\gspace$ and we
can identify the space $\Gamma$ with the right coset space
$G_{0}/G$ for appropriate closed subgroup $G_{0}$ of $G$. We
then build the induced covariant representation following
Mackey's construction (\cite{T}; Theorem 4.2). Our key result
in this section is Theorem ~\ref{thm:inducedrep} regarding
covariant factor representations of $\cstar$-dynamical system
with compact groups. This theorem extends a similar result in
the context of finite groups obtained by Arias and Latremoliere
(\cite{Ar}; Theorem 3.4). As a corollary of Theorem
~\ref{thm:inducedrep} we show that every irreducible
representation of $\cp$ is induced from a stability group.

In Section 5, we study covariant irreducible representation
$(\pi, U)$ of $(A, G_{P},\sigma)$ such that ker($\pi$) = $P$,
where $P\in$ Prim $A$ and show that in this case the
representation $\pi$ of $A$ must be homogeneous. As a
corollary, we get that $\cp$ satisfies the strong-EHI property.

I would like to thank David Pitts and Dana Williams for their
helpful comments regarding this paper. I would especially like
to thank my adviser Allan Donsig for his support and guidance
throughout my graduate studies including the writing of this
paper.

\section{Background}
Suppose that $G$ is a topological (resp. Borel) group; that is,
$G$ is a topological (resp. Borel) space and a group such that
the map $(s, t)\in G\times G\mapsto s^{-1}t \in G$ is
continuous (resp. Borelian). When $G$ is a topological group,
$G$ is often considered as a Borel group equipped with the
Borel structure determined by its topology. Let $\Gamma$ be a
topological (resp. Borel) space. Suppose that an anti
homomorphism of $G$ into the group of all homeomorphisms (resp.
Borel-automorphisms) of $\Gamma$ is given, denoting the
homeomorphism (resp. Borel-automorphism) of $\Gamma$
corresponding to $s\in G$ by $\gamma\in\Gamma \mapsto
\gamma\cdot s\in \Gamma$. If the map: $(\gamma, s)\in
\Gamma\times G\mapsto\gamma\cdot s\in \Gamma$ is continuous
(resp Borelian), then $\Gamma$ is said to be a topological
(resp. Borel) $G$-space. By a measure $\mu$ on a Borel space
$\Gamma$, we shall mean a complete measure determined by a
$\sigma$-finite measure on the Borel sets of $\Gamma$. For each
$s\in G$ define a measure $s(\mu)$ on $\Gamma$ by $s(\mu)
(E)=\mu(E\cdot s)$. We say that $\mu$ is quasi-invariant if
$s(\mu)$ is equaivalent to $\mu$ for each $s\in G$ and we call
the measure space $(\Gamma, \mu)$ a $G$-measure space.

If a quasi-invariant measure $\mu$ on a Borel $G$-space $\Gamma$
satisfies the condition that $\mu (E)=0$ or $\mu (\Gamma - E)=0$ for
every Borel set $E$ of $\Gamma$ with $\mu (E\triangle (E\cdot s))=0$
for every $s\in G$, then $\mu$ is said to be ergodic. Given a unital
$\cstar$-dynamical system $(A, G, \sigma)$ we say the group action
is ergodic if the only fixed elements of $A$ under the group action
are the scalars. Similarly, given a $W^*$-dynamical system $(\bA, G,
\tau)$ (defined at the beginning of section 3) we say that action of
$G$ on the von Neumann algebra $\bA$ is ergodic if the only fixed
elements are the scalars. Note that if $(\Gamma, \mu)$ is a standard
Borel $G$-measure space such that the corresponding action of $G$ on
$\Linf\gspace$ given by $(\tau_{s}f)(\gamma)=f(\gamma\cdot s)$ is
continuous in the strong operator topology (SOT) we can form a
$W^*$-dynamical system $(\Linf\gspace, G, \tau)$. In this case,
$(\Gamma, \mu)$ is an ergodic $G$-measure space if and only if the
action of $G$ on $\Linf\gspace$ is ergodic.

\begin{remark}\label{rem:SOTcts}
The continuity on the measure space can be expressed by the
requirement that $\mu(E\bigtriangleup (E\cdot s))\rightarrow 0$ as
$s\rightarrow e$ for each measurable set $E$ with $\mu(E)< \infty$
(\cite{B}; p.285).
\end{remark}

Let $\gspace$ be a Borel $G$-measure  space. For each
$\gamma\in\Gamma$ define $O_{\gamma}=\{{\gamma\cdot s: s \in G}\}$
to be the orbit of $\gamma$ under the group action. If there is
$\gamma \in \Gamma$ such that $\mu (\Gamma - O_{\gamma})=0$ then
$\gspace$ is said to be transitive. Clearly transitivity implies
ergodicity. As mentioned in the introduction the converse is not
true in general.

\section{Ergodic Actions of Compact Groups}
Let $G$ be a locally compact group,  $\bA$ a von Neumann algebra,
and $\tau$ a point-wise SOT-continuous homomorphism of $G$ into the
automorphism group of $\bA$ then we call the triple $(\bA, G, \tau)$
a $W^*$-dynamical system. It is well known that given a
$W^*$-dynamical system $(\bA, G, \tau)$ the set $\bA^{c}$ of
$x\in\bA$ such that the function $s\mapsto \tau_{s}x$ is norm
continuous is a $G$ invariant $\cstar$-algebra and it is
$\sigma$-weakly dense in $\bA$ (\cite{B}; Proposition III.3.2.4).
Since $\bA^{c}$ is unital it follows from the Double Commutant
Theorem that $\bA^{c}$ is SOT-dense in $\bA$. Using the proof of the
Spectral Theorem we will show that $\bA^{c}$ is equivalent to a
space of continuous functions on a second countable compact
Hausdorff space.

\begin{lemma}\label{lem:SpectralThm}
Let $N$ be a masa on a Hilbert space $\sK$ and $\zeta\in\sK$ be a
cyclic separating vector for $N$. Suppose $M$ is a unital
$\cstar$-subalgebra of $N$ such that $\overline{M\zeta}=\sK$. Then
there exists a compact Hausdorff space $Y$ and a finite Borel
measure $\nu$ and a unitary $V:\sK\rightarrow L^{2} (Y, \nu)$ such
that $VNV^*=\Linf (Y, \nu)$ and $VMV^*=C(Y)$.
\end{lemma}

\begin{proof}
Let $\rho:M\rightarrow C(Y)$ be the Gelfand isomorphism. Define a
positive linear functional $\phi$ on $M$ by $\phi(x)=\langle
x\zeta,\zeta \rangle$. Then there is a finite positive Borel measure
$\nu$ on $Y$ such that
\[ \phi(x)=\int_{Y} \rho(x)d\nu \] for all $x\in M$.

Let $\pi_{\phi}:M\rightarrow \sB (L^{2} (Y, \nu))$ be the
corresponding GNS representation with $\bone$ as the cyclic vector.
Since $\zeta$ is a separating vector then the map
$V:M\zeta\rightarrow \pi_{\phi}(M)\bone$ given by
$V(x\zeta)=\pi_{\phi}(x)\bone$ is well defined. Clearly, $V$ is an
isometry. Hence we can extend $V$ to a unitary from $\sK$ onto
$L^{2} (Y,\nu)$. Moreover, $\pi_{\phi}(x)=VxV^*$ for all $x\in M$ so
that $VMV^*=\pi_{\phi}(M)=C(Y)$. To see that $VNV^*=\Linf (Y, \nu)$
let $x_1 \in M$ and $x_{2}\in N$ then
\[(Vx_{1}V^*)(Vx_{2}V^*)=(Vx_{2}V^*)(Vx_{1}V^*). \]
So $(Vx_{2}V^*)\subseteq(VMV^*)'=(C(Y))'=\Linf (Y,\nu)$. Conversely,
if $T\in\Linf (Y,\nu)\subseteq (VNV^*)'$ then $T(VxV^*)=(VxV^*)T$,
for all $x\in N$. So $x(V^* TV)=(V^* TV)x$, for all $x\in N$. Thus
$V^* TV\in N'=N$ and $T=V(V^* TV)V^* \in VNV^*$.

\end{proof}

Let $\gspace$ be a standard Borel $G$-measure space where the group
action is continuous in the appropriate sense. Consider the
corresponding $W^*$-dynamical system $(\Linf\gspace, G, \tau)$. Then
$\Linf\gspace^{c}$ is SOT-dense in $\Linf\gspace$. Let $\zeta\in
L^2\gspace$ be a cyclic, separating vector for $\Linf\gspace$. Then
we can apply Lemma ~\ref{lem:SpectralThm} to $N=\Linf\gspace$,
$M=\Linf\gspace ^{c}$, and $\zeta$.
\begin{cor}\label{cor:Spectal}
Let $G$ be a locally compact group and let $\gspace$ be a standard
Borel $G$-measure space. Then there is a compact Hausdorff space $Y$
together with a finite positive Borel measure $\nu$ and a unitary
$V: L^2 \gspace \rightarrow L^2 (Y, \nu)$ such that $V \Linf\gspace
V^* = \Linf (Y, \nu)$ and $V \Linf\gspace ^{c} V^* = C(Y)$.
\end{cor}

Consider the $W^*$-dynamical system $(\LY, G, \tau')$ where
$\tau_{s}'(VfV^*)=V(\tau_{s}f)V^*$ for all $s\in G$ and $f\in
\Linf\gspace$. Then by construction we get that $\LY^{c}=C(Y)$. In
particular, $(C(Y),G, \tau')$ is a $\cstar$-dynamical system. Hence
there is an action of $G$ on $Y$ so that $(Y, G)$ is a topological
$G$-space and
\[(\tau_{s}'f)(y)=f(y\cdot s) \]
for all $y\in Y, s\in G$ and $f\in C(Y)$ (\cite{W1}; Proposition
2.7). We would like to show that the above equality holds for all
functions in $\LY$.
\begin{lemma}\label{lem:extend}
In the above situation, let $s\in G$ and $g\in \LY$. Then
$(\tau_{s}'g)(y)=g(y\cdot s)$ for almost all $y\in Y$.
\end{lemma}

\begin{proof}
Let $g$ be in the unit ball of $\LY$ then by the Kaplansky Density
Theorem there is a sequence $(f_{i})$ in the unit ball of $C(Y)$
such that $f_{i}\rightarrow g$ in $\sigma$-SOT. Since $\nu$ is
finite, a simple computation shows that there is a subsequence
$f_{i_{j}}$ converging to $g$ almost everywhere. So without loss of
generality we can assume that $f_{i}\rightarrow g$ almost
everywhere. In particular,\[ f_{i} (y\cdot s)\rightarrow g(y\cdot s)
\] for almost all $y\in Y$. Since automorphisms of von Neumann
algebras are $\sigma$-strong continuous then $f_{i}\rightarrow g$ in
$\sigma$-SOT implies $\tau_{s}'f_{i}\rightarrow \tau_{s}'g$ in
$\sigma$-SOT. By the same argument as above,
\[ (\tau_{s}'f_{i}) (y)\rightarrow (\tau_{s}'g)(y) \] for almost all $y\in Y$. Since $f_{i} \subseteq C(Y)$ then
$(\tau_{s}'f_{i})(y)=f_{i} (y\cdot s)$ for all $y\in Y$ and $i$. It
follows $(\tau_{s}'g)(y)=g(y\cdot s)$ for almost all $y\in Y$.
\end{proof}

\begin{cor}\label{cor:quasiinv}
Let $(Y, \nu)$ be as in Lemma ~\ref{lem:extend}. Then $\nu$ is a
quasi-invariant measure.
\end{cor}
\begin{proof}
Let $X$ be a Borel subset of $Y$. Then $\nu (X)=0 \iff \chi_{X} =0
\iff \tau_{s}'(\chi_{X}) =0 \iff \chi_{(X\cdot s)} =0 \iff \nu
(X\cdot s)=0$
\end{proof}

Suppose the action of $G$ on $\Linf\gspace$ is ergodic then the
action of $G$ on $C(Y)$ must also be ergodic. In general, as
mentioned in the introduction, ergodic actions are far from being
transitive. However, if $G$ is a compact group the two notions
coincide. To this end, we need the following fact which was
initially proved by Albeverio  (\cite{Al}; Lemma 2.1) but we we
offer a different proof.
\begin{lemma}\label{lem:transitive}
Let $G$ be a compact group. Let $X$ be a compact, Hausdorff
topological $G$-space. Suppose the action of $G$ on $C(X)$
given by $(\sigma_{s}f)(x)=f(x\cdot s)$ is ergodic, i.e. the
only $G$ invariant functions are the constant functions. Then
the action of $G$ on $X$ is transitive. Moreover, there exists
a closed subgroup $G_0$ of $G$ such that the right coset space
$G_0/G$ with the quotient topology is homeomorphic to $X$.
\end{lemma}
\begin{proof}
For each $x\in X$ define the orbit of $x$ to be
$O_{x}=\{{x\cdot s: s \in G}\}$. Since the map $s\mapsto x\cdot
s$ is continuous from $G\rightarrow X$ and $G$ is compact then
$O_x$ is compact for each $x\in X$. In particular, $O_x$ is
closed for each $x\in X$.

Fix $x_{0} \in X$. Suppose there is $x_{1} \in X - O_{x_{0}}$ then
$O_{x_{0}}$ and $O_{x_{1}}$ are disjoint closed subsets of $X$. By
Urysohn's Lemma there exists a continuous function $f: X\rightarrow
[0,1]$ such that $f(x_{0}\cdot s)=0$ for all $s\in G$ and
$f(x_{1}\cdot s)=1$ for all $s\in G$. Define a function
$g:X\rightarrow [0,1]$ by $g(x)=\int_{G} f(x\cdot s) dm(s)$. We want
to show that $g$ is continuous. To this end, let $\epsilon
>0$ be given; extend $f$ to $\overline{f}:X\times G \rightarrow
[0,1]$ by defining $\overline{f}(x,s)= f(x\cdot s)$. Then
$\overline{f}$ is continuous function with compact support so we can
find a finite open cover $\{F_i \times G_i\}_{i=1}^n$ of $X\times G$
such that $|f(x\cdot s)-f(y\cdot t)|< \epsilon$ whenever $(x,s)$ and
$(y,t)$ are both in $F_i \times G_i$ for some $i=1,..,n$. Given any
$x\in X$ define $F_x = \bigcap \{F_i : x\in F_i\}$. It is not hard
to check that $|f(x\cdot s)-f(y\cdot s)|< \epsilon$ for all $y\in
F_x$ and $s\in G$. Then $|g(x)-g(y)|\leq \int_{G} |f(x\cdot s) -
f(y\cdot s)| dm(s)  \leq \epsilon$ for all $y\in F_x$. It follows
that $g$ is continuous.

Moreover, $g$ is $G$-invariant and hence must be constant on
$X$. But $g(x_{0})=0$ and $g(x_{1})=1$, contradiction. It
follows that $O_{x_{0}}=X$.

To prove the second part of the statement let $G_{x_0}=\{s\in G: x_0
\cdot s=x_0\}$. Then $G_{x_0}$ is a closed subgroup of $G$ and the
right coset space $G_{x_0}/G$ is compact in the quotient topology.
Moreover, it is easy to see that the map $G_{x_0}\cdot s\mapsto
x_0\cdot s$ is a continuous bijection from $G_{x_0}/G$ onto $X$.
Since $G_{x_0}/G$ is compact and $X$ is Hausdorff it follows that
$G_{x_0}/G$ is in fact homeomorphic to $X$.
\end{proof}


\begin{cor}\label{cor_secondctble}
Let $G$ be a second countable compact group. Let $X$ be a
compact, Hausdorff topological $G$-space. Suppose the action of
$G$ on $C(X)$ given by $(\sigma_{s}f)(x)=f(x\cdot s)$ is
ergodic. Then $X$ is a second countable topological space.
\end{cor}

Applying Lemma~\ref{lem:transitive} to $(C(Y), G, \tau')$ we see
that $G$ acts transitively on $Y$. We are now ready to prove our
main result.
\begin{thm}\label{thm:main}
Let $G$ be a second countable, compact group. Let $\gspace$ be a
standard Borel $G$-measure space. Suppose the action of $G$ on
$\gspace$ is ergodic and the corresponding action of $G$ on
$\Linf\gspace$ is SOT-continuous. Then $G$ acts transitively on
$\gspace$.
\end{thm}
\begin{proof}
We know by Lemma ~\ref{lem:SpectralThm} that there is a
compact, Hausdorff space $Y$ together with a probability
measure $\nu$ and a unitary $V: L^2 \gspace \rightarrow L^2 (Y,
\nu)$ such that $V \Linf\gspace V^* = \Linf (Y, \nu)$ and $V
\Linf\gspace ^{c} V^* = C(Y)$. We define the action of $G$ on
$\LY$ as in Lemma~\ref{lem:extend} then $(Y, \nu)$ becomes a
Borel $G$-measure space by Corollary~\ref{cor:quasiinv}. Since
$G$ is a second countable, compact group then $Y$ is a second
countable topological space by Corollary
~\ref{cor_secondctble}. In particular, $(Y, \nu)$ is a standard
Borel $G$-measure space.

It follows from Mackey's Theorem 5 in \cite{M} that there are
invariant Borel subsets $Y' \subseteq Y$ and $\Gamma ' \subseteq
\Gamma$ and a Borel isomorphism $\theta :Y'\rightarrow \Gamma '$
such that
\begin{enumerate}
\item $\mu(\Gamma - \Gamma')=\nu (Y-Y')=0 $.
\item $\theta(y\cdot s)=\theta(y)\cdot s$ for all $y\in Y', s\in G$.
\end{enumerate}
To show that $\Gamma '$ is an orbit of $G$  let $\gamma_{1},
\gamma_{2} \in \Gamma '$. Let $y_{1}, y_{2} \in Y'$ such that
$\theta(y_{i})=\gamma_{i}$. By Lemma~\ref{lem:transitive} we know
that $G$ acts transitively on $Y'$ so there is $s\in G$ such that
$y_{1}\cdot s=y_{2}$. It follows
$\gamma_{2}=\theta(y_{2})=\theta(y_{1}\cdot s)=\theta(y_{1})\cdot
s=\gamma_{1}\cdot s$ which completes the proof.
\end{proof}

\section{Covariant Representations of Separable $\cstar$-algebras
and Their Compact Automorphism Groups}

In this section we will assume $\cp$ is a separable system and all
Hilbert spaces are separable.  A covariant representation of $\cp$
on a Hilbert space $\sH$ is a pair $(\pi, U)$ where $\pi$ is a
non-degenerate representation of $A$ on $\sH$ and $U$ is a
SOT-continuous homomorphism of $G$ into the unitary group of $\sB
(\sH)$ such that
\[U(s)\pi(a)U(s)^*=\pi(\sigma_{s}a)\]
for all $a\in A$ and $s\in G$.

Let $G_{0}$ be a closed subgroup of $G$ and denote $G_{0}/G$ to be
the corresponding right coset space endowed with the quotient
topology. Let $(\pi_{0}, U_{0})$ be a covariant representation of
$(A, G_{0}, \sigma)$ on a separable Hilbert space $\sHo$. Then
following Mackey's construction of induced representations we can
construct a new covariant representation $(\pi, U)$ of $\cp$, which
is called the induced covariant representation.

In general, if $G$ is a locally compact group then $G_{0}/G$ does
not always admit a $G$-invariant measure so  the construction of
induced representations for groups involves the use of a
quasi-invariant measure $\mu$ on $G_{0}/G$. However, if $G$ is a
compact group there exists a unique, up to scalar multiple,
$G$-invariant Radon measure on $G_{0}/G$ (\cite{F}; Corollary 2.51).
Since the induced representation is independent, up to unitary
equivalence, of the choice of the quasi-invariant measure
(\cite{M1}; Theorem 2.1) the construction of the induced
representation is considerably simplified.

We now describe induced covariant representations following the
construction given in \cite{T}. Let $G_{0}$ be a closed
subgroup of a compact group $G$ and let $(\pi_{0}, U_{0})$ be a
covariant representation of $(A, G_{0}, \sigma)$ on a separable
Hilbert space $\sHo$. Let $\mu$ be a fixed $G$-invariant
measure on $G_{0}/G$. Let $\sH$ denote the induced
representation space then $\sH$ is the space of all $\sHo$
valued functions $\xi$ on $G$ satisfying the following
conditions:
\begin{enumerate}
\item $\langle \xi(s), h_{0}\rangle $ is Borel function of $s$ for all $h_{0}\in \sHo$.
\item $\xi(ts)=U_{0}(t)\xi(s)$ for all $t\in G_{0}$ and all $s\in G.$
\item $\int_{G_{0}/G} \langle \xi(s), \xi(s)\rangle d\mu(\overline{s})<\infty$ .

\end{enumerate}
Define $U$ to be the homomorphism of $G$ into the unitary group of
$\sB (\sH)$ given by:
\[(U(t)\xi)(s)=\xi(st)\] for all
$\xi\in\sH$ and $s, t\in G$. And for each $a\in A$ define an
operator $\pi(a)$ on $\sH$ by:
\[(\pi(a)\xi)(s)=\pi_{0}(\sigma_{s}a)\xi(s)\]
for all $\xi\in\sH$ and $s\in G$. Then $\piU$ is easily checked to
be a covariant representation of $\cp$:
\begin{align*}
U(t)\pi(a)U(t^{-1})\xi(s) & =(\pi(a)U(t^{-1})\xi)(st)\\
&=  \pi_{0}(\sigma_{st}a)(U(t^{-1})\xi)(st)\\
&= \pi_{0}(\sigma_{st}a)\xi(s)= \pi(\sigma_{t}a)\xi(s)
\end{align*} for all $s,t\in G$ and $a\in A$.
Since the $G$-invariant measure $\mu$ is unique up to a scalar
multiple (\cite{F}; Theorem 2.49) the induced representation is
independent of the choice of the measure.

Let $\piU$ be a covariant representation of $\cp$ on $\sH$. We say
that $\piU$ is irreducible if the only operators that commute with
$\pi(a)$ and $U(s)$ for all $a\in A, s\in G$ are the scalars.

Following \cite{T} we define a system of imprimitivity for $\piU$ to
be a commutative von Neumann algebra $\bA$ acting on $\sH$ such
that:
\begin{enumerate}
\item $\bA \subseteq \pi(A)'$.
\item $U(s) \bA U(s)^*= \bA$ for all $s \in G$.
\end{enumerate}
Note that (2) implies that $G$ acts by automorphisms on $\bA$.
Moreover, since $U$ is assumed to be strongly continuous, then for
each $x \in \bA$ the map $s \mapsto U(s) x U(s)^*$ is continuous in
the strong operator topology. If the only $G$ invariant elements of
$\bA$ are scalars then $\bA$ is called an ergodic system of
imprimitivity. In particular, if $\piU$ is an irreducible covariant
representation then $\bA$ is always an ergodic system of
imprimitivity. Given a system of imprimitivity $\bA$ for $\piU$, not
necessarily ergodic, there exists a standard Borel  $G$-measure
space $\gspace$ and an isomorphism $i$ of the algebra $\Linf
\gspace$ onto $\bA$ such that
\[ U(s)i(f)U(s)^*=i(\tau_{s}f) \]
for each $f\in \Linf \gspace$ and $s\in G$ where
$(\tau_{s}f)(\gamma)=f(\gamma\cdot s^{-1})$ (\cite{M}; Theorem 4).
In the above situation we say that the system of imprimitivity $\bA$
for $\piU$ is based on the $G$-measure space $\gspace$ with respect
to $i$. Note that the ergodicity of the system of imprimitivity
$\bA$ is equivalent to that of the action of $G$ on $\gspace$
(\cite{M}; Theorem 3). As in \cite{T} we say that a system of
imprimitivity $\bA$ is transitive if the corresponding Borel
$G$-measure  space is transitive. It follows from Theorem 5 in
\cite{M} that the definition of transitivity is independent of the
choice of $G$-space $\gspace$. Moreover, if a system of
imprimitivity for a covariant representation is transitive then by
(\cite{M2}; Theorem 6.1), the associated $G$-measure space $\gspace$
can be identified with the right coset space $G_{0}/G$ of a closed
subgroup of $G$ together with a $G$-invariant measure on $G_{0}/G$.

If $\bA$ is an ergodic system of imprimitivity for $\piU$  on a
Hilbert space $\sH$ then we can assume
$\sH=L^{2}\gspace\bigotimes \sHo$ and $\bA=\Linf \gspace
\bigotimes I_{\sHo}$ (\cite{M2}; Theorem 5.2). Moreover, the
action of $\bA$ on $\sH$ is given by
\[(i(f)\xi)(\gamma)=f(\gamma)\xi(\gamma)\]
for all $f\in \Linf \gspace$ and $\xi\in L^{2}(\sHo, \Gamma, \mu)$.
In addition,  there exists a Rep$(A: \sHo)$-valued measurable
function $\gamma\in \Gamma \mapsto\pi_{\gamma}\in$ Rep$(A:\sHo)$
such that
\[(\pi(a)\xi)(\gamma)=\pi_{\gamma}(a)\xi(\gamma)\]
for each $a\in A$, $\xi\in\sH$ and almost all $\gamma\in\Gamma$.
Since the action $G$ on $\bA$ is continuous in the strong operator
topology of $\sB(\sH)$ then the corresponding action of $G$ on
$\Linf\gspace$ is also continuous in the strong operator topology of
$\sB (L^{2} \gspace)$. Using Theorem ~\ref{thm:main} we obtain the
following result.
\begin{cor}\label{thm:transsyst}
Let $(\pi, U)$ be an irreducible covariant representation of a
separable system $\cp$ where $G$ is compact. Suppose $\bA$ is a
system of imprimitivity for $\piU$ then $\bA$ is transitive.
\end{cor}

A natural choice for a system of imprimitivity for $\piU$ is the
center of the commutant of $\pi(A)$, which we denote by $\cent$. In
particular, if $\piU$ is a factor representation then $\cent$ is an
ergodic system of imprimitivity for $\piU$. In this case, $\piU$ is
particularly easy to describe. Combining Theorem ~\ref{thm:main} in
Section 3 and Theorem 5.2 in \cite{T} we obtain the following
result.
\begin{thm}\label{thm:inducedrep}
Let $\piU$ be a factor (resp. irreducible) representation of a
separable system $\cp$ where $G$ is compact. Then there exists a
unique closed subgroup $G_{0}$ of $G$ and a unique covariant
representation $(\pi _{0}, U_{0})$ of the subsystem $(A,
G_{0},\sigma)$ such that $\piU$ is induced by $(\pi _{0}, U_{0})$,
where the uniqueness is up to equivalence. Moreover,
\begin{enumerate}
\item $(\pi _{0}, U_{0})$ is a factor (resp. irreducible) representation.
\item $\pi _{0}$ is a factor representation.
\item There is an isomorphism $i:\Linf(G_{0}/G, \mu)\rightarrow \mathbf{Z}(\pi(A)')$ given by
$(i(f)\xi)(s)=f(\overline{s})\xi(s)$.
\end{enumerate}

\end{thm}
Let $G$ be a finite group and $\piU$ be an irreducible
representation of $\cp$. Then we know by the above theorem that
$\piU$ is induced from an irreducible representation $(\pi _{0},
U_{0})$ of $(A, G_{0},\sigma)$ where $\pi_{0}$ is a factor
representation. Define an action of $G_{0}$ on the commutant of
$\pi_{0}(A)$ by $\tau_{s}(T)=U_{0}(s)TU_{0}(s)^*$ for all $s\in
G_{0}$ and $T\in \pi_{0}(A)'$. Since $G$ is finite and acts
ergodically on $\pi_{0}(A)'$ then $\pi_{0}(A)'$ must be finite
dimensional. It follows that $\pi_{0}$ is a direct sum of finitely
many equivalent irreducible representations. Consequently, Theorem
~\ref{thm:inducedrep} can be viewed as a generalization of a similar
result for finite groups obtained by Arias and Latremoliere
(\cite{Ar}; Theorem 3.4).

Let $P$ be a primitive ideal of $A$ and define $G_{P}:= \{s\in G :
\sigma_{s}P=P\}$. Note that $G_{P}$ is a closed subgroup of $G$.
Applying Theorem ~\ref{thm:inducedrep} we get the following
corollary.

\begin{cor}\label{grs}
Let $\piU$ be an irreducible representation of $\cp$. Then there
exists a primitive ideal $P$ of $A$ and a covariant representation
$(\pi _{P}, U_{P})$ of the subsystem $(A, G_{P}, \sigma)$ such that
$\piU$ is induced by $(\pi _{P}, U_{P})$. Moreover, ker $\pi _{P}$ =
$P$.
\end{cor}
\begin{proof}
By Theorem ~\ref{thm:inducedrep} there exists a closed subgroup
$G_{0}$ of $G$ and a unique covariant representation $(\pi _{0},
U_{0})$ of the subsystem $(A, G_{0}, \sigma)$ such that $\piU$ is
induced by $(\pi _{0}, U_{0})$. Since $A$ is separable and $\pi_{0}$
is a factor representation ker$(\pi _{0})$ $\in$ Prim $A$. Let $P$
:= ker $\pi _{0}$ then $G_{0}\subseteq G_{P}$. We take $(\pi _{P},
U_{P})$ to be the representation of $(A, G_{P}, \sigma)$ induced by
the representation $(\pi _{0}, U_{0})$ of the subsystem $(A, G_{0},
\sigma)$.

Moreover, it follows from Lemma ~\ref{lem:homog} in the next section
that ker $\pi _{P}$ = $\bigcap_{r\in G_{P}} \sigma_{r}$(ker $\pi
_{0}$) = $\bigcap_{r\in G_{P}} \sigma_{r}P$ = $P$.
\end{proof}

We note that the above corollary generalizes the GRS Theorem
\cite{G} in the case of compact groups.

\section{Strong EHI}
In this section we continue working with a separable system $\cp$
where $G$ is a compact group. Let $\pi$ be a representation of $A$
on a separable Hilbert space $\sH$. If $E$ is a projection in the
commutant $\pi(A)'$ of $\pi$ then we denote $\pi^{E}$ to be the
subrepresentation of $\pi$ acting on $E\sH$.

Let $G_{0}$ be a closed subgroup of $G$ and $(\pi _{0}, U_{0})$ be a
covariant representation of $(A, G_{0}, \sigma)$ on $\sHo$. Let
$\piU$ be the covariant representation of $\cp$ on $\sH$ induced by
$(\pi _{0}, U_{0})$ then there is a natural family of projections in
$\pi(A)'$ associated with Borel subsets of $G_{0}/G$. Consider the
map $i:\Linf(G_{0}/G, \mu)\rightarrow \pi(A)'$ given by
$(i(f)\xi)(s)=f(\overline{s})\xi(s)$. For each nonzero Borel subset
$E$ of $G_{0}/G$ we denote $\pi^{E}$ to be the subrepresentation of
$\pi$ acting on $i(\chi_{E})\sH$.

\begin{lemma}\label{lem:homog}
In the above situation, let $Q$ := ker $\pi _{0}$. If $F$ is an open
subset of $G_{0}/G$ then ker $\pi^{F} = \bigcap_{s\in q^{-1}(F)}
\sigma_{s^{-1}}Q$.
\end{lemma}

\begin{proof}

Note that the quotient map $q:G\rightarrow G_{0}/G$ is continuous
and open. Let $F$ be an open subset of $G_{0}/G$ and suppose there
is an $a\in A$ such that $a\notin \bigcap_{s\in q^{-1}(F)}
\sigma_{s^{-1}}Q$ then we will show that $\pi^{F}(a)\neq 0$. To this
end, let $s\in q^{-1}(F)$ such that $\pi_{0}(\sigma_{s}a)\neq 0$.
Choose a unit vector $h\in \sHo$ and $\epsilon > 0$ so that
\[\|\pi_{0}(\sigma_{s}a)h\|\geq 2\epsilon\] Then as in the proof of Lemma
6.19 in \cite{W1} we can construct a function $\xi\in C_{b}(G,
\sHo)\subseteq \sH$ such that
\[\|\xi(s)-h\|\leq \epsilon/\|a\|\]
It follows that
$\|\pi_{0}(\sigma_{s}a)\xi(s)-\pi_{0}(\sigma_{s}a)h\|\leq\|\pi_{0}(\sigma_{s}a)\|\cdot\|\xi(s)-h\|\leq\|a\|\cdot(\epsilon/\|a\|)=\epsilon$.
By the reverse triangle inequality we get
\[\|\pi_{0}(\sigma_{s}a)\xi(s)\|\geq\epsilon\]
Since $\pi_{0}(\sigma_{s_{j}}a)\rightarrow \pi_{0}(\sigma_{s}a)$
whenever $s_{j}\rightarrow s$ and $\xi\in C_{b}(G, \sHo)$ there
exists an open neighborhood $F_{s}\subseteq G_{0}/G$ of $G_{0}s$
such that
\[\|\pi_{0}(ta)\xi(t)\|>\epsilon/2\]
for all $t\in q^{-1}(F_{s})$. Then
$\pi^{F}(a)(\chi_{q^{-1}(F_{s}\cap F)}\xi)\neq 0$.

Clearly, $\bigcap_{s\in q^{-1}(F)} \sigma_{s^{-1}}Q\subseteq
\pi^{F}$.
\end{proof}

We call $\pi$ a homogeneous representation if ker $\pi^{E}$ = ker
$\pi$ for every nonzero projection $E\in \pi(A)'$. It follows from
Lemma G.3 in \cite{W1} that $\pi$ is a homogeneous representation if
ker $\pi^{E}$ = ker $\pi$ for every nonzero projection $E\in
\pi(A)'\cap \pi(A)''$. A structure theory developed by Effros in
\cite{E} allows us to decompose arbitrary representations into a
direct integral of homogeneous representations that has very useful
properties. In particular, the following result is due to Echterhoff
and Williams \cite{W2}:
\begin{thm}\label{thm:EW}
Let $\cp$ be a separable system. Suppose that $\rho$ is a
homogeneous representation of $A$ with ker $\rho$ = P, and that
$\rho\times_{\sigma} V$ is an irreducible representation of
$A\times_{\sigma} G_{P}$. Then the representation of
$A\times_{\sigma} G$ induced by $\rho\times_{\sigma} V$ is
irreducible.
\end{thm}

We would like to use Theorem ~\ref{thm:EW} to prove the strong-EHI
property for separable systems involving compact groups. To this end
we prove the following theorem.
\begin{thm}\label{thm:homog}
Let $\cp$ be a separable system where $G$ is a compact group.
Suppose $P$ is a primitive ideal of $A$ and $\piU$ is an irreducible
covariant representation of $(A, G_{P}, \sigma)$ such that ker
$\pi$= $P$. Then $\pi$ is a homogeneous representation of $A$.
\end{thm}

\begin{proof}
Note that $G_{P}$ is a closed subgroup of $G$ so $G_{P}$ is compact.
We know by Theorem ~\ref{thm:inducedrep} that there exists a closed
subgroup $G_{0}$ of $G_{P}$ and an irreducible covariant
representation $(\pi _{0}, U_{0})$ of the subsystem $(A, G_{0})$
such that $\piU$ is induced by $(\pi _{0}, U_{0})$. Moreover, there
is an isomorphism $i:\Linf(G_{0}/G_{P}, \mu)\rightarrow
\mathbf{Z}(\pi(A)')$ given by $(i(f)\xi)(s)=f(\overline{s})\xi(s)$.
Let $E$ be a Borel subset of $G_{0}/G_{P}$ of nonzero measure. Let
$\sH$ be the representation space of $\piU$ and denote $\pi^{E}$ to
be the subrepresentation of $\pi$ acting on $i(\chi_{E})\sH$. Then
it is enough to show that ker $\pi^{E}$ = ker $\pi$ .

Denote $Q$ := ker $\pi _{0}$.  If $F$ is an open subset of
$G_{0}/G_{P}$ denote $F':=\{s^{-1} : s\in q^{-1}(F)\}$. Then by
Lemma ~\ref{lem:homog}, ker $\pi^{F} = \bigcap_{s\in F'}
\sigma_{s}Q$. Since $G_{0}/G_{P}$ is compact there is
$\{t_{j}\}_{1\leq j \leq n} \subseteq G_{P}$ such that
$G_{P}=\bigcup t_{j} F'$. Then $P$ = ker $\pi$ = $\bigcap_{r\in
G_{P}} \sigma_{r}Q$ = $\bigcap \sigma_{t_{j}} (\bigcap_{s\in F'}
\sigma_{s}Q)$ = $\bigcap \sigma_{t_{j}}$ (ker $\pi^{F}$). Since $P$
is a prime ideal and $P$ is $G_{P}$-invariant it follows that $P$ =
ker $\pi^{F}$. In particular, $\|\pi^{F}(a)\|=\|\pi(a)\|$ for all
$a\in A$.

Now let $K$ be a compact subset of $G_{0}/G_{P}$ of nonzero measure.
By a simple compactness argument we can find $G_{0}s\in K$ such that
every open neighborhood of $G_{0}s$ intersects with $K$ in a set of
positive measure. Then by the arguments similar to Lemma
~\ref{lem:homog} it follows that ker $\pi^{K}\subseteq$ ker
$\pi_{0}\circ s$. We want to show that ker $\pi_{0}\circ s
\subseteq$ ker $\pi$. To this end, suppose $\pi_{0}(\sigma_{s}a)=0$
and let $\epsilon > 0$ be given. Since
$\pi_{0}(\sigma_{s_{j}}a)\rightarrow 0$ whenever $s_{j}\rightarrow
s$ we can find an open neighborhood $F'$ of $s$ in $G_{P}$ such that
$\|\pi_{0}(\sigma_{t}a)\|<\epsilon$ for all $t\in F'$. Then
$\|\pi(a)\|=\|\pi^{q(F')}(a)\|<\epsilon$. Thus $\pi(a)=0$ as
claimed. It follows ker $\pi^{K}= P$.

Finally, if $E$ a nonzero Borel subset of $G_{0}/G_{P}$ then we can
choose a compact subset $K\subseteq E$ such that $\mu(K)> 0$.
Suppose $\pi^{E}(a)=0$ then $\pi^{K}(a)=0$. It follows
$\|\pi(a)\|=\|\pi^{K}(a)\|=0$. So ker $\pi^{E}=P$.
\end{proof}

Combining Theorem ~\ref{thm:EW} and Theorem ~\ref{thm:homog} we
obtain the following important corollary.

\begin{cor}\label{cor:ehi}
Let $\cp$ be a separable $\cstar$-dynamical system where $G$ is
compact. Then $\cp$ satisfies the strong-EHI property.
\end{cor}

As mentioned in the introduction it remains unknown whether the
strong-EHI property holds for an arbitrary $\cstar$-dynamical
system. We can inquire about a weaker property of $\cstar$-dynamical
systems, called simply the EHI property, where we ask every
primitive ideal of $A\times_{\sigma} G$ to be induced from a
stability group (see \cite{W2}). However, even with an additional
assumption that $G$ is amenable it is not known whether all
separable $\cstar$-dynamical systems satisfy the EHI property.

\end{document}